\documentclass{article}

\usepackage{arxiv}
\usepackage[utf8]{inputenc} 
\usepackage[T1]{fontenc}    
\usepackage{hyperref}       
\usepackage{url}            
\usepackage{booktabs}       
\usepackage{amsfonts}       
\usepackage{nicefrac}       
\usepackage{microtype}      
\usepackage{graphicx}
\usepackage{subcaption}
\usepackage{doi}
\usepackage{amsmath}
\usepackage[dvipsnames]{xcolor}
\usepackage{todonotes}
\usepackage[table]{xcolor}
\usepackage{caption}
\usepackage{siunitx}
\usepackage{comment}

\usepackage[
backend=biber,
style=numeric-comp,
sorting=none,
giveninits=true
]{biblatex}
\def\liquid{\ell}



\title{Heat and mass transfer through fabric: a model for fabric drying with heated cylinders}

\author{ {S. Bellavia\thanks{Full Professor of Numerical Analysis}} \\
	Department of Industrial Engineering\\
	University of Florence\\
	Firenze, Viale Morgagni 40-44\\
	\texttt{stefania.bellavia@unifi.it} \\
    \And
    {N. Fiorini\thanks{Pre-doc}} \\
	Department of Industrial Engineering\\
	University of Florence\\
	Firenze, Viale Morgagni 40-44\\
	\texttt{nicolo.fiorini@unifi.it} \\
    \And
    {A. Milazzo \thanks{Associate Professor of Thermodynamics}} \\
	Department of Industrial Engineering\\
	University of Florence\\
	Firenze, Via di S. Marta 3\\
	\texttt{adriano.milazzo@unifi.it} \\
    \And
    {A. Papini\thanks{Associate Professor of Numerical Analysis}} \\
	Department of Industrial Engineering\\
	University of Florence\\
	Firenze, Viale Morgagni 40-44\\
	\texttt{alessandra.papini@unifi.it} \\
}

\begin{document}

\maketitle

\begin{abstract}
Textile drying is a key operation in the textile production cycle as it represents one
of the most energy-intensive stages and plays a critical role in determining both product quality and overall process
efficiency.
In this work we propose a mathematical model  for the drying process of a generic textile material using  heated cylinders, operating under low-pressure conditions. The model's parameters are estimated by nonlinear least squares regression.
Given a specific fabric, the developed model allows to predict the drying time and the residual moisture content.
The model is validated using real world data provided by a major  Italian textile company. 
\end{abstract}


\section{Introduction}
Within the textile industry, the reduction of both the extensive consumption of water resources and carbon footprint  associated with the  fabric washing process constitutes a multifaceted challenge, that necessitates a fundamental rethinking of the processes involved. 

As highlighted in the review  \cite{DEFRAEYE2014323}, the improvement of energy efficiency in thermal processes—such as drying—plays a crucial role in this transition. In particular, the development of predictive and physics-based models is increasingly recognized as a key enabling factor for the optimization of existing technologies and for the design of novel solutions aimed at waste heat recovery, reduction of processing times, and minimization of product losses, thereby lowering the embodied energy of textile products.

In this perspective, the availability of reliable models capable of predicting drying dynamics represents a valuable tool for industrial process optimization. Accurate predictions of drying times and residual moisture content allow for a more efficient scheduling of production cycles, avoidance of over-drying, and reduction of unnecessary energy consumption. Ultimately, such predictive capabilities support informed control strategies that can significantly contribute to energy savings and to the overall sustainability of textile finishing operations.

In this context, a collaborative effort has been established and strengthened among \textit{S. Stefano S.p.A.}—a distinguished company in the Italian textile sector—, \textit{Lafer S.p.A.}—a leading manufacturer of industrial textile machinery—and the Department of Industrial Engineering at the University of Florence. The objective of this partnership is the development of an industrial-scale system for the washing and drying of fabrics that entirely eliminates the use of water and detergents, substituting them with carbon dioxide (\(\mathrm{CO_2}\)) as the sole processing medium.

Within this broader framework, the present work focuses on the drying phase of the process, which plays a critical role in determining both product quality and overall process efficiency. As a first step in this direction, we present a mathematical model for the drying of a generic textile material using heated cylinders operating under low-pressure conditions. The model is designed to capture the coupled evolution of temperature and moisture content within the fabric thickness, accounting for phase-change effects and for the dependence of thermophysical properties of the wet fabric, such as specific heat capacity, density and thermal conductivity, with respect to the local moisture level.

The model includes a limited number of physically meaningful parameters, namely the evaporation rate coefficient, the residual (non-evaporable) moisture content, and the parameters governing the smooth activation of evaporation with respect to temperature and moisture thresholds. These parameters are estimated through bound-constrained nonlinear least-squares regression using real-world experimental data provided by the industrial partners. The identification procedure explicitly incorporates prior physical knowledge of the system by imposing suitable bounds on the admissible parameter values. The resulting nonlinear least-squares problem is solved using a trust-region algorithm \cite{BM,trust,MMP}.

Since the aforementioned machine is still at the prototype stage, experimental data from this specific system are currently unavailable. Nevertheless, the mathematical model accounts for the different thermo-physical properties of a generic fluid permeating the fabric and has been validated using experimental data collected from a similar drying machine, in which the fabric is washed with water and subsequently dried by means of heated cylinders. 
Once experimental data from the CO$_2$-based machine become available, the model can be adapted by replacing the thermo-physical properties of water with those of carbon dioxide, as well as by updating the corresponding boundary and initial conditions.

The model we have devised in this paper is inspired by already existing models for paper drying \cite{cepitis, pang1997relationship}.
We elaborated on them taking into account the different thermo-physical properties of the material, the specific operating conditions, and the geometry of the system. The model can be used to predict the drying time of a given texture and the average residual moisture content at the end of the drying process carried out by a machine that employs heated rollers.

\subsection{Related works}

The drying process of a given material has been investigated in several studies.
Early foundational work \cite{fourt1951rate} established the baseline kinetics of fabric drying, demonstrating that drying rates are governed primarily by air resistance and fabric thickness rather than the specific type of fiber.  
Haghi et al. \cite{haghi2002mechanism} detailed how moisture diffuses differently in highly hygroscopic fibers compared to weakly hygroscopic ones, while Ismail et al. \cite{ismail1988heat} developed mathematical models to predict the effective thermal conductivity of woven fabrics.  The issue of predicting the effective thermal conductivity and thermal resistance of the woven fabric has also been faced \cite{siddiqui2013finite} by using finite element methods.  
 
Building on these thermodynamic considerations, recent research has mainly focused on the simulation of convective drying processes. Comprehensive mathematical models  \cite{johann2014mathematical,dornyak2025mathematical,sousa2006analysis} to analyze moisture distribution, local thermophysical parameters, and the impact of operational variables during textile drying under convective air flow have been proposed. 
In particular, the models presented by Johan et al. \cite{johann2014mathematical} and Dornyak et al. \cite{dornyak2025mathematical}, describe the convective drying process of textile materials and allow the estimation of local thermophysical and structural parameters of the fabric, while Sousa et al.\cite{sousa2006analysis} investigate experimentally the drying process of crude cotton fabric, by exposing textile samples to a controlled hot air flow within a drying chamber. 
These modeling approaches have also been extended to consumer applications, such as the simulation of domestic drum dryers\cite{wei2017mathematical}.

The heat and mass transfer proces with a single-cylinder drying machine has also been investigated \cite{tran2022parametric} for seven different fabric types. 
More generally, the mathematical principles governing heat and mass transfer in fibrous porous media are shared across different industrial processes. For example, Pang et al. \cite{pang1997relationship} investigated the relationship between transport and diffusion models in softwood drying, while Ghosh et al. \cite{ghosh2011fundamentals} and Cepitis \cite{cepitis} discussed the theoretical foundations and mathematical modeling of industrial paper drying systems.

 Despite the extensive literature on textile drying, most existing models describe convective drying processes in which air is used as the primary heat and mass transfer medium. The drying process considered in the present work differs significantly from this configuration. In the proposed system, the fabric is dried through direct contact with heated cylinders operating in a low-pressure environment, a mechanism that is conceptually closer to industrial paper drying processes, but here applied to textile materials.
Another distinguishing aspect of the proposed model concerns the estimation of the thermophysical properties of the wet fabric. In contrast with several approaches in the literature, where the solid matrix and the liquid phase are treated separately, the present model considers the wet fabric as an effective single system composed of fibers and liquid. Consequently, the thermodynamic parameters used in the model represent effective properties of the combined fabric--fluid system.

\subsection{Paper's organization}
The paper is organized as follows.
Section \ref{problem description} provides a detailed description of the industrial drying process and of the physical system under consideration, together with the modeling assumptions adopted in this work. In Section \ref{mathmodel}, the mathematical model is derived, including the coupled energy and mass balance equations, the constitutive relations for the material properties, and the associated initial and boundary conditions, which together define the resulting initial–boundary value problem. Section \ref{numresults} is devoted to the numerical assessment of the model and to the  parameter identification process, where the model's parameters are estimated by fitting  experimental drying data through a nonlinear least-squares approach. Finally, conclusions and possible directions for future work are discussed in the last section.

\section{Problem description}
\label{problem description}
\begin{figure}
    \centering
    \includegraphics[width=0.7\linewidth]{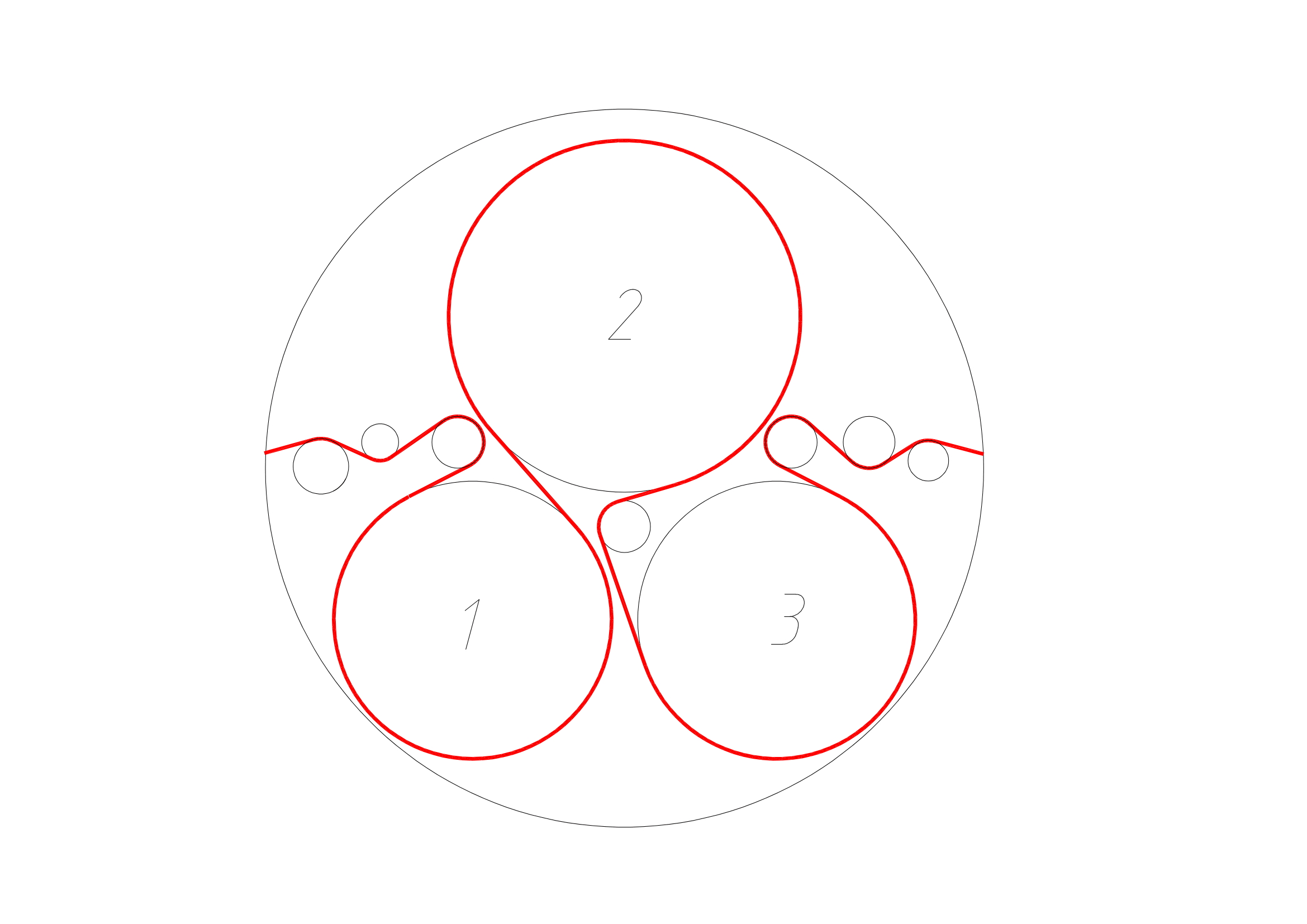}
    
    \caption{Cross-sectional view of the machine. The fabric is highlighted in red, while the rollers with a larger diameter correspond to the heated cylinders.
}
\label{fig:sezione_rullo}
\end{figure}
In order to derive the mathematical model for fabric drying, the process through which the fabric is introduced into the machine and subsequently dried must first be described. The main components involved in the drying stage are therefore outlined below; a cross-sectional view of the system is reported in Figure~\ref{fig:sezione_rullo}.

The system consists of a vacuum chamber that maintains a constant pressure of \(0.1\,\text{bar}\). Under these operating conditions, the evaporation temperature of the fluid permeating the fabric is significantly reduced. Inside the chamber, three heated rollers are employed for the drying process, together with seven additional feed rollers. The latter ensure the continuous motion of the fabric within the machine and maintain the required tension, allowing the fabric to remain in close contact with the heated cylinders at all times.

To improve the efficiency of the drying process, the fabric is alternately exposed to the lower and upper heated rollers, so that both sides of the material are progressively dried. The heated rollers are maintained at a constant temperature \(T_{\mathrm{cyl}}\) by the machine control system, while the fabric translation speed within the chamber is also kept constant.

Based on these technical considerations, the drying process is modeled by focusing on the interaction between the fabric and a \textbf{single heated cylinder}, with the aim of deriving the corresponding system of partial differential equations. Although the industrial machine is equipped with three heated cylinders, the drying dynamics associated with each cylinder can be described using the same mathematical formulation.
Once the temporal evolution of the relevant quantities has been determined for a single cylinder, the solution obtained at the exit of the first cylinder can be used as the initial condition for the second one, and the same procedure can be repeated for the third cylinder. In this way, the overall drying process inside the machine is reconstructed through a sequential application of the single-cylinder model.

The quantities we want to model in time and space are the \textit{temperature} of the wet fabric \(T(x,t)\) and the \textit{moisture concentration} \(M(x,t)\), a dimensionless variable defined as the ratio between the local mass of liquid \(m_\liquid(x,t)\)  and the constant mass of fibers \(m_f\):
\begin{equation*}
    M(x,t) := \frac{m_\liquid(x,t)}{m_f}.
\end{equation*}

After defining the relevant quantities, the physical model is formulated under a set of simplifying assumptions, which allow the problem to be reduced to a tractable framework.

It is assumed that the heated rollers do not experience significant temperature variations due to contact with the wet fabric, as the machine control system compensates for the associated heat losses and maintains the cylinders at a constant temperature.

Furthermore, resistance to vapor transport through the fabric is neglected. Once the liquid evaporates, the generated vapor is assumed to escape freely into the surrounding environment. In addition, liquid motion within the fabric due to capillary or gravitational effects is not considered. These assumptions are supported by the relatively short duration of the drying phase and the small thickness of the fabric, which limit the relevance of internal transport phenomena.

Assuming that the fabric is isotropic, the spatial domain of the problem can be restricted to the fabric thickness only. The model can therefore be formulated as a one-dimensional problem in the spatial variable \(x \in (0, L)\), measured along the direction normal to the surface of the heated cylinder.

Under these assumptions, the drying process can be described by a system of two coupled equations, namely the \textit{energy balance equation} and the \textit{mass balance equation}. The proposed formulation is valid for a generic liquid permeating the fabric. In the experimental setup considered in this study, however, the fabrics are wetted with water. Consequently, the simulations, the experimental analysis and the specific initial and boundary conditions reported in the following sections are carried out using the thermodynamic properties of water.

\section{The mathematical model}
\label{mathmodel}
The present section illustrates the mathematical model describing the drying process of the fabric.

The model describes the behavior of two quantities: the fabric temperature and the moisture content within the fabric, when exposed to a constant heat source provided by the heated cylinder. The system’s evolution is formulated as an initial boundary value problem, expressed through a set of two coupled partial differential equations. Estimates of the physical and thermodynamic parameters of the wet fabric, obtained during the study, are incorporated within these equations.

\subsection{Energy balance equation: heat transfer}
The heat from the cylinder surface through the fabric is described by the following heat equation: 

\begin{equation*}
            c\rho\frac{\partial T}{\partial t} = \frac{\partial}{\partial x}\left[ \lambda\frac{\partial T}{\partial x} \right] + h_{\liquid,v}\rho_\liquid\frac{\partial M}{\partial t},
\end{equation*}

a non-homogeneous parabolic partial differential equation, where $h_{\liquid,v}$ is the latent heat of vaporization for the liquid, $\rho_\liquid$ is the density of liquid and the sink term \( h_{\liquid, v }\rho_\liquid  \frac{\partial M(x,t)}{\partial t} \) accounts for the energy required to evaporate the mass of liquid.

The terms \(c,\, \lambda,\, \rho \), which respectively denote the specific heat capacity, the thermal conductivity and the density of the wet fabric, are not constant, since they all depend on the liquid content within the fabric. 
As there are no explicit expressions for these quantities, the following estimates are adopted:

\begin{itemize}
    \item \textbf{Heat capacity.} Weighted average of the heat capacities of the liquid, \(c_\liquid\), and the fabric, \(c_f\), with respect to the corresponding masses \(m_\liquid,m_f\):
    \begin{equation*}
        c(M) := \frac{m_\liquid c_\liquid + m_f c_f}{m_\liquid + m_f} = \frac{Mc_\liquid + c_f}{M+1}
    \end{equation*}
    \item \textbf{Thermal conductivity.} The estimate given in Slavinec et al. \cite{thermalcond} is used; the quantity increases linearly with \(M\) for every \(M<M_c\), where $M_c$ is a critical moisture level experimentally evaluated\footnote{For values of $M$ greater than the critical level $M_c$, it is observed that the thermal conductivity of the fabric is equal to the thermal conductivity of the liquid itself.} at \(M_c = 2\):
    \begin{equation*}
        \label{eq:thermal_conductivity}
            \lambda(M) := M\frac{\lambda_\liquid - \lambda_f}{M_c} + \lambda_f, \quad 0 < M \leq M_c
    \end{equation*}
where $\lambda_\ell,\lambda_f$ are the thermal capacities of the liquid and the fabric, respectively;
    \item \textbf{Density.} The average density of the fabric is expressed as
    \begin{equation*}
            \label{eq:density}
            \rho(M) = \frac{m_\liquid + m_f}{V_\liquid +V_f} = 
            \frac{(M+1)\rho_\liquid\rho_f}{M\rho_f + \rho_\liquid}
        \end{equation*}
        where $V_\liquid, V_f$ are the volumes of the liquid and the fabric, respectively.
\end{itemize}

\subsection{Mass balance equation: moisture evolution}
The equation describing the evolution of moisture within the fabric must satisfy several physical constraints arising from the nature of the process. 

It is assumed that the liquid contained in the fabric remains still during the process, thereby neglecting both capillary and gravitational effects. Mathematically, this assumption translates into the absence of spatial derivatives in the governing equation, since in every point of the considered domain the liquid flux is zero.

Moreover, a mass loss due to evaporation occurs; since \( m_f \) represents a constant (positive) quantity, the time derivative of the moisture content \(M\) must therefore be non-positive, that is,

\[
\frac{\partial M}{\partial t} \leq 0.
\]

In addition to these phenomena, two important properties of water evaporation must be taken into account: firstly, to completely dry a piece of fabric (\(M = 0\)), 
both extremely high temperatures and prolonged drying times would be required, conditions beyond the scope of the considered industrial drying machine. Secondly, even if the fabric were fully dried, exposure to standard environmental conditions\footnote{20°C, 50\% relative humidity in the air} would cause some water reabsorption from the surrounding environment; so this behavior has to be incorporated into the model.

In addition to that, the evaporation process is assumed to occur only when the temperature of the wet fabric \(T\) reaches the evaporation temperature of the liquid \( T_{\mathrm{evap}} \). Therefore, the partial differential equation describing the mass balance under the aforementioned assumptions takes the form

\begin{equation*}
    \frac{\partial M}{\partial t} =  -\kappa(M) \max(T-T_{evap}, 0) 
\end{equation*}

where \(\kappa: \mathbb{R}^2 \rightarrow \mathbb{R}^+ \) is a non negative function estimating the \textbf{evaporation rate}.

Since the \(\max\) function is not differentiable everywhere, a smooth \(C^{\infty}\) approximation is introduced through logistic-type function:

\begin{equation*}
\label{lnT}
    \max(T - T_{\mathrm{evap}}, 0) \;\approx\; 
    \frac{T - T_{\mathrm{evap}}}{1 + e^{-\beta (T - T_{\mathrm{evap}})}}
\end{equation*}

where \(\beta\) is a positive parameter governing the steepness of the logistic transition. Its value must be consistent with both the physical characteristics of the system and the numerical stability of the solver, in particular with respect to overflow constraints.
The value of \(\beta\) adopted in the simulations is \(\beta = 3\). This choice ensures a suitable approximation of the maximum function while preventing numerical overflow. 
Indeed, the largest negative temperature deviation occurs at the beginning of the drying process:
\begin{equation*}
    \label{T_bound}
\min_{(x,t)}\!\big(T(x,t) - T_{\mathrm{evap}}\big)
   = T(x,0) - T_{\mathrm{evap}}
\end{equation*}
and at $t = 0$  it is desirable that
\begin{equation*}
\label{approximation}
\frac{T(x,0) - T_{\mathrm{evap}}}
     {1 + e^{-\beta (T(x,0) - T_{\mathrm{evap}})}} 
     \approx 0.
\end{equation*}
In our simulations $T(x,0)=T_0$ is constant and $T_0 - T_{\mathrm{evap}}=-30$, so that 
$$
\frac{T_0 - T_{\mathrm{evap}}}{1 + e^{-\beta (T_0 - T_{\mathrm{evap}})}} 
\simeq -2.5\times 10^{-38} ~~~~{\rm and}~~~~ e^{-\beta (T_0 - T_{\mathrm{evap}})} \simeq 1.2 \times 10^{39},
$$
safely far from the IEEE floating-point double-precision overflow threshold.

A more accurate calibration of this parameter could be achieved by  
estimating  $\beta$
 through
a nonlinear least-squares approach as for the other parameters of the model 
(see \eqref{problem_ls}-\eqref{problem_ls_bound} below).
 On the other hand, the experimental results obtained with the selected value are satisfactory, and including \(\beta\) in the parameter space would significantly increase the computational cost of the parameter estimation procedure.

In a similar fashion, another logistic like function is chosen to estimate correctly the form for $\kappa(M)$; it is shown \cite{cepitis} that during the drying process there is an initial phase of constant drying, then the intensity of the process slows down as we reach the minimum residual moisture content. So the following form for $\kappa(M)$ is chosen

\begin{equation*}
    \kappa (M) := \frac{k}{1 + e^{- \gamma(M-M_b)}}
\end{equation*}
with
\begin{itemize}
    \item \(k\), representing the evaporation rate; larger values of \(k\) correspond to a more intense drying process;
    \item \(M_b\), representing the water level at which the drying process has to slow down;
    \item \(\gamma\), representing the change in intensity during the evaporation process as we approach the value $M_b$ for $M$; larger values of $\gamma$ correspond to a more constant drying intensity.
\end{itemize}

As described in the following section, the parameters $k, M_b, \gamma$ are   estimated by fitting experimental data.

To conclude, the final form of the mass balance equation can be written as

\begin{equation*}
    \frac{\partial M}{\partial t} = - \frac{k}{1 + e^{- \gamma(M-M_b)}} \frac{T-T_{evap}}{1 + e^{-\beta(T-T_{evap})}} 
\end{equation*}

\subsection{Initial Boundary Value Problem}
To solve the two coupled partial differential equations, suitable initial and boundary conditions are required. The initial temperature \(T_0\) and moisture \(M_0\) of the fabric can be measured before the drying process begins:
\begin{equation*}
    \label{initial conditions}
    T_0 = T(x,0), \quad M_0 = M(x,0),\quad 0<x<L.
\end{equation*}

Concerning the boundary conditions, it can be assumed that the temperature at the contact interface between the heated cylinder and the fabric remains constant and equal to that of the cylinder itself. Moreover, due to the substantial temperature difference between the cylinder and the fabric, 
and considering that the water evaporation temperature is reduced to approximately \(45\,^{\circ}\mathrm{C}\) under low-pressure conditions \footnote{In our specific case the pressure is approximately $0.1$ bar}, an instantaneous vaporization at the same contact interface may be reasonably assumed.
 
This leads to the following Dirichlet boundary conditions:
\begin{equation*}
    \label{DIRboundary}
    T(0,t) = T_{cyl}, \quad M(0,t) \approx 0, \quad t>0.
\end{equation*}

On the opposite side of the fabric, namely the surface not in contact with the heated cylinder, a convective heat exchange with the surrounding environment must be considered. This gives rise to the following Robin boundary condition:
\begin{equation*}
\label{ROBboundary}
    \frac{\partial T}{\partial x} = -\frac{z_{ht}}{\lambda(M)} (T(L,t) - T_{env}) 
\end{equation*}

where \(z_{ht} \in \mathbb{R}^+\) is a heat transfer coefficient.

The $z_{ht}$ heat transfer coefficient has been evaluated combining convection and radiation 
as $z_{ht} = h_{conv} + h_{irr}$. 
The convection coefficient $h_{conv}$ has been calculated as 
\[
h_{conv}=\frac{Nu \cdot \lambda_{vap} }{D}
\] 

where $ \lambda_{vap} = 0.02$ is the thermal conductivity of water vapour, $D$ is the cylinder diameter and the Nusselt number $Nu$  is given\cite{kreith2010principles} by
\[ 
Nu = 0.11\,(\,0.5 \,Re^2 + Gr \cdot Pr\,)^{0.35}
\]

where $Re$, $Pr$ and $Gr$ are the numbers of Reynolds, Prandtl and Grashow respectively.

The radiation coefficient $h_{irr}$ has been obtained  assuming a simplified geometry, where the cylinder is surrounded by a cylindrical enclosure. In this way we may use a view-factor $F_{env \rightarrow cyl}\,$ equal to the ratio between the diameters of the cylinder and the enclosure. Referring to $q$, the heat flux per unit area of the cylinder surface, one has: 

\begin{equation*}
q = \frac{ \sigma (T^4_{cyl} - T^4_{env}) }{ \frac{1 - \epsilon_{cyl}}{\epsilon_{cyl}} + 1 + \frac{F_{env \rightarrow cyl}(1 - \epsilon_{env})}{\epsilon_{env}} }
\end{equation*}

where $\sigma$ is the Stefan-Boltzmann constant and $\epsilon$ is the surface emissivity.
Then 
\begin{equation*}
    h_{irr} = 
    \frac{q}{T_{cyl} - T_{env}} = 
    \frac{ \sigma ( T_{cyl}^2 + T_{env}^2 )( T_{cyl} + T_{env} ) }{ \frac{1 - \epsilon_{cyl}}{\epsilon_{cyl}} + 1 + \frac{ F_{env \rightarrow cyl} (1-\epsilon_{env})}{\epsilon_{env}}}.
\end{equation*}

Given these conditions, the complete Initial and Boundary Value Problem can be stated as

\begin{equation}
\begin{cases}
     \dfrac{\partial T}{\partial t} = \dfrac{1}{c(M)\rho(M)}\left[ \dfrac{\partial}{\partial x}\left[ \lambda(M)\dfrac{\partial T}{\partial x} \right] + h_{\liquid,v}\rho_\liquid\dfrac{\partial M}{\partial t}\right],  &x\in(0,L),\,t>0\\[1.4em]
    \dfrac{\partial M}{\partial t} =  - \dfrac{k}{1 + e^{- \gamma(M-M_b)}} \dfrac{T-T_{evap}}{1 + e^{-\beta(T-T_{evap})}} , &x\in(0,L),t>0\\[1.4em]
    T(x,0) = T_0, \,\,\, 
    M(x,0) = M_0, & \text{Initial conditions}\\[.75em]
    T(0,t) = T_{cyl}, \,\,\,
    M(0,t) \approx 0, & \text{Boundary conditions}   \\[.75em]
    \dfrac{\partial T}{\partial x}\Big|_{(L,t)} = - \dfrac{z_{ht}}{\lambda(M)}(T(L,t)-T_{env}), & \text{Flux conditions}  
\end{cases}
\label{eq:IBVP}
\end{equation}

\section{Experimental results}
\label{numresults}
In this section, the experimental results and the methodologies employed in the study are presented.

To solve the initial–boundary value problem \eqref{eq:IBVP}, it is necessary to estimate the model parameters $k$, $M_b$ and $\gamma$ on the basis of the available experimental data.
The dataset consists of pairs \(\{\theta_i,  M_{\tau,i}\}_{i = 1 \dots N}\); each element of the dataset corresponds to a single fabric sample and $N$ is the total number of distinct samples.
In particular 
\begin{equation*}
    \label{eq:dataset pairs}
    (\theta_i , M_{\tau,i}) \equiv ( (\tau_i, L_i, T_{cyl,i}, M_{0,i}), M_{\tau,i})
    \in \mathbb{R}^4 \times \mathbb{R} , \quad i = 1\dots N,
\end{equation*}
where, for each distinct sample, \(\tau\) denotes the total duration of the drying process, \(L\) is the thickness of the fabric, \(T_{cyl}\) is the temperature of the three heated cylinders,  \(M_0\) is the initial moisture content and \(M_\tau\) denotes the
final average moisture content of the fabric:
\begin{equation*}
   M_{\tau,i} = \frac{1}{L_i} \int_0^{L_i} M_i(x,\tau_i) dx,
\end{equation*}
where $M_i(x,t)$ is the moisture concentration of the i-th sample. Note that this distribution is unknown and 
the experimental measures for $M_{\tau,i}$ were made using only the total weight of the wet and dried samples.



We decided to choose the free parameters by a least-squared minimization of the residual between the model and the data, $f : \mathbb{R}^3 \rightarrow \mathbb{R}^{N}$, defined as

\begin{equation*}
    \label{f_ls}
    f(k, M_b, \gamma):= ( M_{\tau,1} - \hat{ M}_{\tau,1},\dots,  M_{\tau,N} - \hat{ M}_{\tau,N}),
\end{equation*}

where $\hat{M}_{\tau, i}$ is the prediction for $M_{\tau,i}$ computed by solving the IBVP \eqref{eq:IBVP} with the choice of parameters $k, M_b, \gamma$, i.e.
\begin{equation*}
    \hat{M}_{\tau, i} = \frac{1}{D} \sum_{j = 1}^D \hat{M}_i(x_j, \tau), 
\end{equation*}
$\{ x_j \}_{j = 1, \dots, D}$ is an equally-spaced grid  of the interval $[0,L]$ 
and $\{\hat{M}_i(x_j, \tau)\}_{j = 1, \dots, D}$ is the numerical solution of the IBVP for the fabric sample $i$. The solution is calculated by a finite element solver implemented in Python using the FEniCS library \cite{fenicsLoggEtal2012, fenicsKirbyLogg2006, fenicsAlnaesEtal2015}.

More precisely,
the optimal parameters for the model are determined by solving
the following bound-constrained nonlinear least-squares  problem:
\begin{gather}
\min_{k, M_b, \gamma} \| f(k, M_b, \gamma) \|_2^2  \label{problem_ls}\\
\text{with } 
k \in [\,10^{-4}, 10^{-3}\,], \quad
M_b\in [\,0.02, 0.2\,], \quad\gamma \in \,[10, 150\,].\label{problem_ls_bound}
\end{gather}

Taking into account that $\kappa, M_b$ and $\gamma$  have physical meaning only if they belong to suitable intervals, 
the box constraints for the parameters have been established by looking for patterns in experimental data; regarding this matter, it is important to note that we are looking for a \textbf{unique} set of "optimal" parameters which are independent from the type of fabric, but in principle these parameters, which govern the dynamics of the drying process, could be fine-tuned for every different type of fabric.
Since such a detailed, fabric-dependent calibration is beyond the scope of the present study, the bounds imposed on such parameters - especially on \(M_b\) - are intentionally chosen to be relatively broad.

In our experimentation, we used the Trust Region Reflective method \cite{trust} implemented in the SciPy module \textit{optimize.least\_square} to solve the bound-constrained least-squares problem \eqref{problem_ls}-\eqref{problem_ls_bound}.
It is important to note that since the gradient of the objective function is unknown, it is necessary to calculate a finite-difference approximation of it at each iteration of the trust region method.
This implies that at iteration $j$ of the method we need to evaluate the function $f$ four times\footnote{We need to evaluate $f$ in  $(k_j, M_{b,j}, \gamma_j), (k_j + h, M_{b,j}, \gamma_j), (k_j, M_{b,j} + h, \gamma_j), (k_j, M_{b,j}, \gamma_j + h)$ where $h$ is the increment for each parameter, which is needed to compute the corresponding component of the gradient through finite-difference approximation.} to compute the approximation of the gradient and the value of the function in the current choice of the parameters $(k_j, M_{b,j}, \gamma_j)$. 
This procedure can be computationally demanding, since we need to solve the IBVP $N$ times for each evaluation of $f$, which adds up to solving the differential problem $4N$ times for each iteration of the trust-region method.

The complete results for all tested samples are reported in Tables \ref{tab:results0}, \ref{tab:results1} and \ref{tab:results2}. The \textbf{KEY} column contains the unique identifier of each sample, while the \textbf{TRUE} and \textbf{PRED} columns report the values of \(M_{\tau,i}\) and \(\hat{M}_{\tau, i}\), respectively.
In the \textbf{UNDER-OVER DRIED} column we can see if the predicted value matches the experimental value: if $|M_{\tau,i} - \hat{M}_{\tau, i}| \leq 0.015$ the prediction on the sample is considered correct and the \textit{correctly dried} tag is attached in the column. Conversely if $M_{\tau,i} - \hat{M}_{\tau, i}> 0.015$ or $M_{\tau,i} - \hat{M}_{\tau, i} < - 0.015$ the \textcolor{red}{over-dried}/\textcolor{blue}{under-dried} tags are respectively attached. 

The results obtained with the 
parameter values provided by  the trust-region optimization procedure are 
reported in Table \ref{tab:results0}. 
The values of such parameters are  
listed in the Table's caption. They have been computing using $ (k_0, M_{b,0}, \gamma_0) = ( 5\times 10^{-4},\,  0.1 , \, 70 ) $ as initial guess for the trust-region procedure.

In Tables \ref{tab:results1} and \ref{tab:results2} we report the results obtained using alternative parameter configurations in order to show that the performance  of the model are not strongly dependent on the parameters' choice. In fact, comparable values of the mean squared error (MSE) and mean absolute error (MAE) were observed, indicating a consistent predictive performance.

\begin{figure}[!h]
\centering
\begin{subfigure}{0.49\linewidth}
    \centering
    \includegraphics[width=\linewidth]{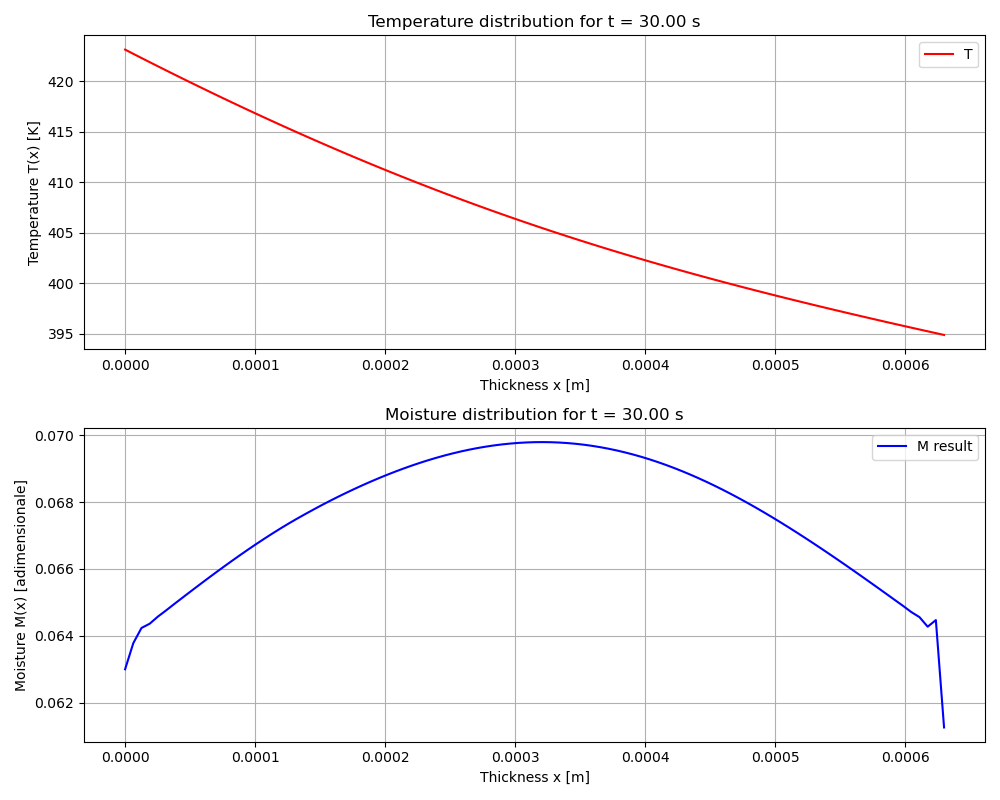}
    \caption{Sample 3}
    \label{fig:sample2_MT_process}
\end{subfigure}
\hfill
\begin{subfigure}{0.49\linewidth}
    \centering
    \includegraphics[width=\linewidth]{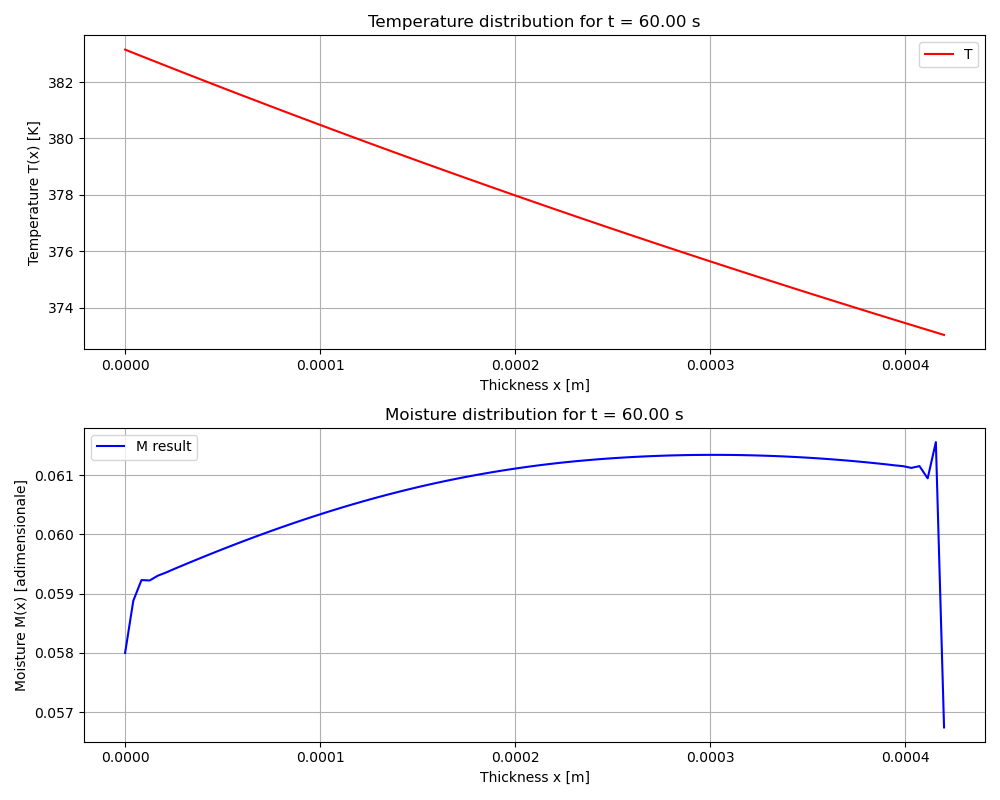}
    \caption{Sample 14}
    \label{fig:sample13_MT_process}
\end{subfigure}

\caption{\textcolor{red}{Temperature} and \textcolor{blue}{Moisture} distributions across the fabric thickness at the end of the drying process for two representative samples.}
\label{fig:MT_process_samples}
\end{figure}

Plots of $T$ and $M$ at the end of the drying process for two specific samples are displayed in Figures \ref{fig:sample2_MT_process} and \ref{fig:sample13_MT_process}.
Figure \ref{fig:sample2_MT_process} shows the temperature and moisture distributions across the fabric thickness at $t = 30s$ for the Sample 3. The temperature decreases almost monotonically from the cylinder side ($x=0$) to the external surface ($x=L$). This behaviour is a direct consequence of the boundary conditions: the temperature is fixed at the cylinder–fabric interface, while the opposite surface exchanges heat with the environment through a convective–radiative condition. Heat therefore propagates through the thickness by conduction, while part of the supplied energy is consumed by water evaporation through the latent heat term in the energy balance. 
The moisture profile shows a bell-shaped distribution, with lower values near both surfaces and a slightly higher value in the interior. This behaviour can be explained by the operating sequence of the drying process. The fabric undergoes three drying phases and is flipped after each phase: it is first dried on side A, then on side B, and finally again on side A. As a consequence, both surfaces are directly exposed to the heated cylinder during the process and experience stronger evaporation, while the inner region remains comparatively less heated.
We can also identify a similar  behaviour for Sample 14 in Figure \ref{fig:sample13_MT_process} where the temperature decreases almost linearly across the thickness and the  moisture distribution curve is flatter, with respect to Sample 3. This flatness is expected and is the result of various concurring factors: this sample was dried for $60s$ instead of $30s$, the thickness is different ($L_3 =  6.3 \times 10^{-4}, \; L_{14} = 4.2 
\times 10^{-4}$ ) and the initial water quantity is almost the same ($M_{0,3} = \num{0.63}, \, M_{0,14} = \num{0.58}$ ), which results in the sample being more uniformly dried.

Overall, the results show that the model captures the general behaviour of the drying process in a consistent and physically meaningful way. Although a limited number of samples are not classified correctly, these discrepancies are most likely attributable to the relatively small amount of available experimental data rather than to systematic deficiencies of the model. In fact, the Mean Absolute Error remains below 0.015 for the best set of parameters, indicating that the predictions are generally accurate and that the adopted parameter set provides a satisfactory global description of the process across different fabric types.

\begin{table}[]
\centering
\rowcolors{2}{white}{gray!5}
\begin{tabular}{cccccc}
\rowcolor{gray!30}
\textbf{KEY} & \textbf{TRUE} & \textbf{PRED} & \textbf{RELATIVE ERROR} & \textbf{ABSOLUTE ERROR} & \textbf{UNDER-OVER DRIED} \\
\midrule
1 & 4,76E-02 & 5,23E-02 & 0,0886 & 0,0046 & correctly dried \\
2 & 4,35E-02 & 5,21E-02 & 0,1648 & 0,0086 & correctly dried \\
3 & 6,84E-02 & 6,54E-02 & 0,0467 & 0,0031 & correctly dried \\
4 & 4,55E-02 & 6,05E-02 & 0,2486 & 0,0150 & correctly dried \\
5 & 5,88E-02 & 5,97E-02 & 0,0153 & 0,0009 & correctly dried \\
6 & 4,60E-02 & 5,95E-02 & 0,2277 & 0,0136 & correctly dried \\
7 & 7,78E-02 & 5,53E-02 & 0,4062 & 0,0225 & \textcolor{red}{over-dried} \\
8 & 5,19E-02 & 5,43E-02 & 0,0456 & 0,0025 & correctly dried \\
9 & 6,25E-02 & 5,28E-02 & 0,1833 & 0,0097 & correctly dried \\
10 & 9,29E-02 & 5,95E-02 & 0,5614 & 0,0334 & \textcolor{red}{over-dried} \\
11 & 4,52E-02 & 5,73E-02 & 0,2105 & 0,0121 & correctly dried \\
12 & 6,21E-02 & 6,55E-02 & 0,0516 & 0,0034 & correctly dried \\
13 & 6,86E-02 & 5,24E-02 & 0,3094 & 0,0162 & \textcolor{red}{over-dried} \\
14 & 7,84E-02 & 5,83E-02 & 0,3448 & 0,0201 & \textcolor{red}{over-dried} \\
15 & 3,79E-02 & 6,62E-02 & 0,4275 & 0,0283 & \textcolor{blue}{under-dried} \\
16 & 5,56E-02 & 5,36E-02 & 0,0361 & 0,0019 & correctly dried \\
17 & 3,33E-02 & 6,68E-02 & 0,5012 & 0,0335 & \textcolor{blue}{under-dried} \\
\bottomrule
\end{tabular}
\caption{Comparison between experimental values (TRUE) and  numerically simulated values (PRED) and the associated error metrics. 
Results were obtained by solving the IBVP with the following choice of parameters $k = 9.99 \times 10 ^{-4}$, $M_b = 9.75 \times 10^{-2}$, and $\gamma = 1.49 \times 10 ^2$.
\textbf{Mean Squared Error} : $0.000289$ , \textbf{Mean Absolute Error} : $0.0135$
}
\label{tab:results0}
\end{table}

\begin{table}[]
\centering
\rowcolors{2}{white}{gray!5}
\begin{tabular}{cccccc}
\rowcolor{gray!30}
\textbf{KEY} & \textbf{TRUE} & \textbf{PRED} & \textbf{RELATIVE ERROR} & \textbf{ABSOLUTE ERROR} & \textbf{UNDER-OVER DRIED} \\
\midrule
1 & 4,76E-02 & 4,32E-02 & 0,1023 & 0,0044 & correctly dried \\
2 & 4,35E-02 & 4,31E-02 & 0,0083 & 0,0004 & correctly dried \\
3 & 6,84E-02 & 6,79E-02 & 0,0084 & 0,0006 & correctly dried \\
4 & 4,55E-02 & 5,67E-02 & 0,1983 & 0,0112 & correctly dried \\
5 & 5,88E-02 & 5,53E-02 & 0,0639 & 0,0035 & correctly dried \\
6 & 4,60E-02 & 5,49E-02 & 0,1618 & 0,0089 & correctly dried \\
7 & 7,78E-02 & 4,83E-02 & 0,6100 & 0,0295 & \textcolor{red}{over-dried} \\
8 & 5,19E-02 & 4,62E-02 & 0,1226 & 0,0057 & correctly dried \\
9 & 6,25E-02 & 4,38E-02 & 0,4263 & 0,0187 & \textcolor{red}{over-dried} \\
10 & 9,29E-02 & 5,48E-02 & 0,6945 & 0,0381 & \textcolor{red}{over-dried} \\
11 & 4,52E-02 & 5,10E-02 & 0,1138 & 0,0058 & correctly dried \\
12 & 6,21E-02 & 6,78E-02 & 0,0838 & 0,0057 & correctly dried \\
13 & 6,86E-02 & 4,32E-02 & 0,5897 & 0,0255 & \textcolor{red}{over-dried} \\
14 & 7,84E-02 & 5,28E-02 & 0,4863 & 0,0257 & \textcolor{red}{over-dried} \\
15 & 3,79E-02 & 6,59E-02 & 0,4252 & 0,0280 & \textcolor{blue}{under-dried} \\
16 & 5,56E-02 & 4,55E-02 & 0,2202 & 0,0100 & correctly dried \\
17 & 3,33E-02 & 6,89E-02 & 0,5163 & 0,0356 & \textcolor{blue}{under-dried} \\
\bottomrule
\end{tabular}
\caption{Comparison between experimental values (TRUE) and  numerically simulated values (PRED) and the associated error metrics. 
Results were obtained by solving the IBVP with the following choice of parameters $k = 5 \times 10 ^{-4}$, $M_b = 1 \times 10^{-1}$, and $\gamma = 1 \times 10 ^2$.
\textbf{Mean Squared Error} : $0.000379$ , \textbf{Mean Absolute Error} : $0.0151$
}
\label{tab:results1}
\end{table}

\begin{table}
\centering
\rowcolors{2}{white}{gray!5}
\begin{tabular}{cccccc}
\rowcolor{gray!30}
\textbf{KEY} & \textbf{TRUE} & \textbf{PRED} & \textbf{RELATIVE ERROR} & \textbf{ABSOLUTE ERROR} & \textbf{UNDER-OVER DRIED} \\
\midrule
1 & 4,76E-02 & 5,47E-02 & 0,1291 & 0,0071 & correctly dried \\
2 & 4,35E-02 & 5,45E-02 & 0,2019 & 0,0110 & correctly dried \\
3 & 6,84E-02 & 6,76E-02 & 0,0127 & 0,0009 & correctly dried \\
4 & 4,55E-02 & 6,29E-02 & 0,2774 & 0,0174 & \textcolor{blue}{under-dried} \\
5 & 5,88E-02 & 6,21E-02 & 0,0534 & 0,0033 & correctly dried \\
6 & 4,60E-02 & 6,19E-02 & 0,2576 & 0,0160 & \textcolor{blue}{under-dried} \\
7 & 7,78E-02 & 5,78E-02 & 0,3468 & 0,0200 & \textcolor{red}{over-dried} \\
8 & 5,19E-02 & 5,67E-02 & 0,0861 & 0,0049 & correctly dried \\
9 & 6,25E-02 & 5,52E-02 & 0,1318 & 0,0073 & correctly dried \\
10 & 9,29E-02 & 6,19E-02 & 0,5006 & 0,0310 & \textcolor{red}{over-dried} \\
11 & 4,52E-02 & 6,79E-02 & 0,3343 & 0,0227 & \textcolor{blue}{under-dried} \\
12 & 6,21E-02 & 5,48E-02 & 0,1337 & 0,0073 & correctly dried \\
13 & 6,86E-02 & 6,07E-02 & 0,1302 & 0,0079 & correctly dried \\
14 & 7,84E-02 & 6,85E-02 & 0,1445 & 0,0099 & correctly dried \\
15 & 3,79E-02 & 5,60E-02 & 0,3241 & 0,0182 & \textcolor{blue}{under-dried} \\
16 & 5,56E-02 & 6,92E-02 & 0,1967 & 0,0136 & correctly dried \\
17 & 3,33E-02 & 7,97E-02 & 0,5817 & 0,0464 & \textcolor{blue}{under-dried} \\
\bottomrule
\end{tabular}
\caption{Comparison between experimental values (TRUE) and  numerically simulated values (PRED) and the associated error metrics. 
Results were obtained by solving the IBVP with the following choice of parameters
$k = 1 \times 10^{-3}$, $M_b = 1 \times 10 ^ {-1}$, and $\gamma = 1.5 \times 10^2$.
\textbf{Mean Squared Error} : $0.000328$ , \textbf{Mean Absolute Error} : $0.0144$
}
\label{tab:results2}
\end{table}

\newpage

\section{Conclusions and future work}
In this work, a mathematical model for the drying process of wet fabrics in a vacuum-assisted industrial machine has been proposed. The formulation is based on a coupled system of energy and mass balance equations describing the evolution of temperature and moisture concentration within the fabric thickness.

The model has been developed for a single heated cylinder, which represents the fundamental drying unit of the machine. This formulation allows the overall process involving multiple cylinders to be reconstructed sequentially by using the solution obtained at each stage as the initial condition for the subsequent one.

Model parameters have been estimated using experimental data collected from different fabric samples. The results show that the proposed formulation is capable of reproducing the observed drying behavior with satisfactory accuracy, while maintaining a relatively low computational cost. In addition, the analysis highlights a certain robustness of the model with respect to variations in the selected parameters.

Despite the simplifying assumptions adopted in the present study, the model captures the main physical mechanisms governing the drying process and provides a useful framework for the analysis and optimization of the machine operating conditions.

Future developments of the present model may focus on relaxing some of the simplifying assumptions introduced in this study, with the aim of achieving a more comprehensive description of the drying process in porous textile materials. 

A first direction concerns the extension of the mass balance to include liquid transport within the fabric, accounting for capillary and gravitational effects. This would naturally lead to the introduction of additional transport terms in the governing equations. In a similar way, a more detailed representation of the vapor phase could be considered, including its generation and motion through the porous structure of the material.

Another relevant aspect involves the mechanical interaction between the fabric and the heated cylinders. The expansion of vapor inside the fabric may locally reduce the contact with the cylinder surface, leading to partial detachment and consequently modifying the dominant heat transfer mechanisms, which could then include radiative and convective contributions in addition to conduction.

From a modeling perspective, the assumption of isotropic material properties could also be relaxed. Considering anisotropic fabrics would require the formulation of two- or three-dimensional models capable of accounting for the internal structure and geometry of the textile.

Finally, with specific reference to the industrial machine analyzed in this work,  different operational stages of the process could be incorporated into the model. In particular, the transitions between successive cylinders may introduce intermediate cooling effects before the fabric reaches the next heated roller. Including these mechanisms would allow the model to better reproduce the actual operating conditions of the system and improve its predictive capability.

\section{Acknowledgments}

The authors would like to acknowledge Simone Zumatri\footnote{M. Sc. Simone Zumatri, R\&D Department Manager, Rifinizione S. Stefano S.p.A., zumatri@santostefano.prato.it} from the R\&D Department of  Rifinizione S. Stefano S.p.A. for his valuable support in sample collection and weighing, as well as for providing detailed information regarding fabric typologies and their properties.

The authors also gratefully acknowledge Denise Mantoan\footnote{Ph.D. Denise Mantoan, R\&D Department, Lafer S.p.A., d.mantoan@laferspa.com} and Giorgia Lera\footnote{M. Sc. Giorgia Lera, R\&D Department, Lafer S.p.A., g.lera@laferspa.com} from the R\&D Department of Lafer S.p.A. for their support concerning machine operation, discussions on issues encountered during the design phase, and insightful exchanges on process dynamics.

Finally, the authors would like to thank Stefano Caporali\footnote{Ph.D. Stefano Caporali, Associate Professor, University of Florence, stefano.caporali@unifi.it}, from the Department of Industrial Engineering of the University of Florence, for carrying out the thickness measurements of the samples.

\newpage
\section*{Nomenclature}

\begin{tabular}{ll}
\multicolumn{2}{l}{\textit{Independent variables and geometry}} \\
$x$ & spatial coordinate through fabric thickness (m) \\
$t$ & time (s) \\
$L$ & fabric thickness (m) \\
$D$ & cylinder diameter (m) \\[4pt]

\multicolumn{2}{l}{\textit{State variables}} \\
$T(x,t)$ & fabric temperature (K) \\
$M(x,t)$ & moisture content (-) \\[4pt]

\multicolumn{2}{l}{\textit{Material properties}} \\
$c(M)$ & specific heat capacity of wet fabric (J\,kg$^{-1}$\,K$^{-1}$) \\
$\lambda(M)$ & thermal conductivity of wet fabric (W\,m$^{-1}$\,K$^{-1}$) \\
$\rho(M)$ & density of wet fabric (kg\,m$^{-3}$) \\
$h_{\ell,v}$ & latent heat of vaporization (J\,kg$^{-1}$) \\[4pt]

\multicolumn{2}{l}{\textit{Process temperatures}} \\
$T_{\mathrm{cyl}}$ & cylinder temperature (K) \\
$T_{\mathrm{env}}$ & environment temperature (K) \\
$T_{\mathrm{evap}}$ & evaporation temperature (K) \\[4pt]

\multicolumn{2}{l}{\textit{Evaporation parameters}} \\
$k$ & evaporation rate coefficient \\
$M_b$ & residual moisture content (-) \\
$\beta$ & logistic temperature parameter (K$^{-1}$) \\
$\gamma$ & logistic moisture parameter (-) \\[4pt]

\multicolumn{2}{l}{\textit{Heat transfer quantities}} \\
$z_{ht}$ & overall heat transfer coefficient (W\,m$^{-2}$\,K$^{-1}$) \\
$h_{\mathrm{conv}}$ & convective heat transfer coefficient (W\,m$^{-2}$\,K$^{-1}$) \\
$h_{\mathrm{irr}}$ & radiative heat transfer coefficient (W\,m$^{-2}$\,K$^{-1}$) \\
$q$ & heat flux (W\,m$^{-2}$) \\
$\sigma$ & Stefan--Boltzmann constant (W\,m$^{-2}$\,K$^{-4}$) \\[4pt]

\multicolumn{2}{l}{\textit{Dimensionless numbers}} \\
$Nu$ & Nusselt number \\
$Re$ & Reynolds number \\
$Pr$ & Prandtl number \\
$Gr$ & Grashof number \\
\end{tabular}

\printbibliography
\end{document}